\documentclass[a4paper, 10pt]{article}
\usepackage[cp1251]{inputenc}
\usepackage[russian]{babel}
\usepackage{mathrsfs}
\usepackage{latexsym}
\usepackage{amssymb}
\usepackage{amsmath}
\usepackage{amsopn}
\usepackage{amsthm}
\usepackage{geometry}

\makeatletter
\renewcommand{\@biblabel}[1]{#1.}
\makeatother

\newtheorem{theorems}{Theorem}[section]
\newtheorem{statement}[theorems]{Proposition}
\newtheorem{consequence}[theorems]{Corollary}

\begin{document}

\begin{center}
\textbf{The complete system of differential invariants of a curve in pseudo-euclidean space}\\
\end{center}

\begin{center}\small
V.I.Chilin\\
Department of Mathematics\\
National University of Uzbekistan,\\
Vuzgorodok, 100174, Tashkent, Uzbekistan\\
e-mail: chilin@ucd.uz\\
\phantom{A}\par
K.K.Muminov\\
Department of Mathematics\\
National University of Uzbekistan,\\
Vuzgorodok, 100174, Tashkent, Uzbekistan\\
e-mail: m.muminov@rambler.ru
\end{center}

%%%%%%%%%%%%%%%%%%%%
%%%%%%%%%%%%%%%%%%%%%%%%%%%%%%%%%%%%%%%%
%%%%%%%%%%%%%%%%%%%%%%%%%%%%%%%%%%%%%%%%%%%%%%%%%%%%%%%%%%%%
%%%%%%%%%%%%%%%%%%%%%%%%%%%%%%%%%%%%%%%%%%%%%%%%%%%%%%%%%%%%%%%%%%%%%%%%%%%%%%%%
\renewcommand{\abstractname}{Abstract}
\begin{abstract}
A complete system of differential invariants for equivalence of
curves in the $n$-dimensional pseudo-euclidean space with respect to the
action of each of the groups $K^n \lhd O(n,p,K)$, $K^n \lhd SO(n,p,K)$,
$O(n,p,K)$, and $SO(n,p,K)$, where $K = \mathbf{R}$, or $K =
\mathbf{C}$, and $O(n,p,K)$ $\bigl($respectively, $SO(n,p,K)\big)$
are a pseudo-orthogonal (special pseudo-orthogonal) groups, is given.
\end{abstract}

\phantom{A}


\phantom{A}


%%%%%%%%%%%%%%%%%%%%
%%%%%%%%%%%%%%%%%%%%%%%%%%%%%%%%%%%%%%%%
%%%%%%%%%%%%%%%%%%%%%%%%%%%%%%%%%%%%%%%%%%%%%%%%%%%%%%%%%%%%
%%%%%%%%%%%%%%%%%%%%%%%%%%%%%%%%%%%%%%%%%%%%%%%%%%%%%%%%%%%%%%%%%%%%%%%%%%%%%%%%
\section{Introduction}
\large
\hspace*{\parindent} Let $X$ be an $n$-dimensional linear space over
the field $K$, where $K$ is either the field of real numbers
$\mathbf{R}$ or the field of complex numbers $\mathbf{C}$. Let
$GL(n,K)$ be the group of all invertible linear mappings of space
$X$, and let $G$ be a subgroup of $GL(n,K)$. Curves $\alpha$ and $\beta$ in $X$
are said to be $G$-equivalent if $g(\alpha) =
\beta$ for some $g \in G$.

Finding necessary and sufficient
conditions for $G$-equivalence of given curves $\alpha$ and
$\beta$ is a classical problem of differential geometry. In recent
years this study involved the use of the
theory of invariants and its methods, with the help of which
the differential fields of all $G$-invariant differential rational
functions over paths were studied and finite rational bases of these
fields were described. Knowing these bases allows one to formulate effective
criteria for $G$-equivalence of paths and curves. This approach
was employed in \cite{liter_1} and \cite{liter_2} to solve the
problem of the equivalence of curves under movements in
$\mathbf{R}^n$ and for actions of semi-direct products
$\mathbf{R}^n \lhd SL(n,\mathbf{R})$ of groups $\mathbf{R}^n$ and
$SL(n,\mathbf{R})$, in particular. In \cite{liter_3} this method was
used to solve the problem of the equivalence of curves in a
specific case of the action of the simplectic group
$Sp(2n,\mathbf{C})$, whereas in \cite{liter_4} it was used to solve this problem for actions of
the groups $\mathbf{R}^n \lhd O(n,\mathbf{R})$ and $\mathbf{R}^n \lhd
SO(n,\mathbf{R})$. Furthermore,  in papers \cite{liter_5} and \cite{liter_6},
necessary and sufficient conditions for the
equivalence of paths under the action of pseudo-orthogonal group
$O(n,p,\mathbf{R})$ and special pseudo-orthogonal group
$SO(n,p,\mathbf{R})$ were found.

The present article is devoted to a study of $G$-equivalence
of curves in the space $X$, where $G$ is one of the groups $K^n \lhd
O(n,p,K)$, $K^n \lhd SO(n,p,K)$, $O(n,p,K)$, or $SO(n,p,K)$.

%%%%%%%%%%%%%%%%%%%%
%%%%%%%%%%%%%%%%%%%%%%%%%%%%%%%%%%%%%%%%
%%%%%%%%%%%%%%%%%%%%%%%%%%%%%%%%%%%%%%%%%%%%%%%%%%%%%%%%%%%%
%%%%%%%%%%%%%%%%%%%%%%%%%%%%%%%%%%%%%%%%%%%%%%%%%%%%%%%%%%%%%%%%%%%%%%%%%%%%%%%%
\section{Preliminaries}

\hspace*{\parindent} Let $X$ be an $n$-dimensional linear space over the field $K$, where $K = \mathbf{R}$ or  $K = \mathbf{C}$. An element from $X$ will be represented as $n$-dimensional column vector $x = \{x_i\}_{i=1}^n$, where $x_i \in K$, $i=\overline{1,n}$. Let $(g,x) \rightarrow gx$ be the left action in $X$ of the group $GL(n,K)$ of all invertible linear transformations of $X$, i.e. usual multiplication of matrix $g$ and column vector $x$.

Futher by $I$ there will be denoted an interval $(a,b)$ from $\mathbf{R}$ (the cases $a = -\infty$ and $b = +\infty$ are possible).

An $I$-path in $X$ is a vector-fucntion $x(t) = \{x_j(t)\}_{j=1}^n$ from $I$ into $X$ such that all coordinate functions $x_j: I \rightarrow K$ are infinitely differentiable functions, $j=\overline{1,n}$.

The vector $\mathop{x}\limits^{(r)}(t) =
\big\{\mathop{x_1}\limits^{(r)}(t), \ldots, \mathop{x_n}\limits^{(r)}(t)\big\}$ is called $r$-th derivative of the $I$-path $x(t) = \{x_1(t), \ldots, x_n(t)\}$, where $\mathop{x_i}\limits^{(r)}(t)$ is an
$r$-th derivative of coordinate functions $x_i(t)$, $i = \overline{1,n}$. It is clear that $\mathop{x}\limits^{(r)}(t)$ is also a $I$-path for all
$r = 1, 2, \ldots$

An $I$-path $x(t)$ is called regular, if $\mathop{x}\limits^{(1)}(t) \not= 0$ for all $t \in I$.

Since $X$ is a finite-dimensional space, in  $X$ exists an unique topology $\tau$ such that  $(X,\tau)$ is a Hausdorff topological vector space.
In particular, we may assume that the topology $\tau$ is induced by Euclidean norm  $\|x\| = \left(\sum\limits_{i=1}^n|x_i|^2\right)^{1/2}$.

An $I$-path $x(t)$ is called simple, if the mapping $x: I \rightarrow X$ is injective and the inverse mapping $x^{-1}:\big(x(I),\tau\big) \rightarrow I$ of $x$ is continuous $\big($\cite{liter_7}, page 22$\big)$. Note, that there exists an $I$-path which is not simple and $x(t)$ is injection $\big($see e.g.
\cite{liter_7}, pages 22-23$\big)$. For simple regular paths the following important theorem holds.

%%%%%%%%%%%%%%%%%%%%%%%%%%%%%%%%%%%%%%%%%%%%%%%%%%
%%%%%%%%%%%%%%%%%%%%%%%%%%%%%%%%%%%%%%%%%%%%%%%%%%%%%%%%%%%%%%%%%%%%%%%%%%%%%%%%%%%%%%%%%%%%%%%%%%%%
\begin{theorems} {\rm$\big($\cite{liter_7}, page 25$\big)$} \label{teor_2_1}
Let $x(t)$ $\big($respectively, $y(t)$ $\big)$ be a simple regular  $I_1$-path (respectively, $I_2$-path). If $x(I_1) = y(I_2)$, then there exists a $C^\infty$-diffeomorphism $\varphi: I_2 \rightarrow I_1$ such that $y(t) = x\big(\varphi(t)\big)$ and $\varphi'(t) \not= 0$ for all $t \in I_2$.
\end{theorems}
%%%%%%%%%%%%%%%%%%%%%%%%%%%%%%%%%%%%%%%%%%%%%%%%%%%%%%%%%%%%%%%%%%%%%%%%%%%%%%%%%%%%%%%%%%%%%%%%%%%%
%%%%%%%%%%%%%%%%%%%%%%%%%%%%%%%%%%%%%%%%%%%%%%%%%%

For every $I$-path $x(t) = \{x_j(t)\}_{j=1}^n$, denote by $M(x)(t)$ the $n \times n$-matrix
$\left(x(t)\mathop{x}\limits^{(1)}(t)\ldots\mathop{x}\limits^{(n-1)}(t)\right)$, where $j$-th column has coordinates $\mathop{x_i}\limits^{(j-1)}(t)$,
$i = \overline{1,n}$. The matrix $\left(\mathop{x}\limits^{(1)}(t)\mathop{x}\limits^{(2)}(t)\ldots\mathop{x}\limits^{(n)}(t)\right)$ is denoted by $M'(x)$. An $I$-path $x(t)$ is said to be strongly regular, if the determinant $\det M(x)(t)$ is nonzero for all $t \in I$. If $x(t)$ is not a regular path, then
$\mathop{x}\limits^{(1)}(t_0) = 0$ for some $t_0 \in I$, therefore $\det M(x)(t_0) = 0$, i.e. the path $x(t)$ is not strongly regular. Hence, every strongly regular $I$-path is regular. Note, that in \cite{liter_1}, \cite{liter_2}, \cite{liter_3}, \cite{liter_5},
\cite{liter_6} the term "regular path" means strongly regular path in the above sense.

Let $G$ be a subgroup of the group $GL(n,K)$. Two $I$-paths $x(t)$ and $y(t)$ are called $G$-equivalent, if there exists an element $g \in G$ such that
$y(t) = gx(t)$ for all $t \in I$. Clearly, in this case, $\mathop{y}\limits^{(j)}(t) = g \mathop{x}\limits^{(j)}(t)$, $j = 1, 2, \ldots$, therefore $M(y)(t) = g M(x)(t)$. It is also obvious that, if the latter equality is valid, then $I$-paths $x(t)$ and $y(t)$ are $G$-equivalent.

Consider two bilinear forms defined on $X = K^n$:
\[ (x, y) = x_1 y_1 + \ldots + x_n y_n, \]
\[ [x, y] = x_1 y_1 + \ldots + x_p y_p - x_{p+1}y_{p+1} - \ldots - x_n y_n, \]
where $p \in \{1, \ldots, n-1\}$. Denote by $E$ the identity of group $GL(n, K)$ and by $E_p = (E_{ij}^{(p)})_{i,j=1}^n$ the matrix from $GL(n, K)$, which satisfies $E_{ii}^{(p)} = 1$ for
 $i = 1, 2, \ldots, p$, $E_{ii}^{(p)} = -1$ for $i = p+1, p+2, \ldots, n$, $E_{ij}^{(p)} = 0$ for $i \not= j$.

\begin{sloppypar}
Let $O(n,K)$ $\big($respectively, $O(n,p,K)$ $\big)$ be an othogonal (respectively, pseudo-orthogonal) subgroup of $GL(n,K)$, i.e. $O(n,K) =
\big\{ g \in GL(n,K): g^T g = E \big\} = \big\{g \in GL(n,K): (gx,gy) = (x,y) \text{ for all } x,y \in X\big\}$ $\big($respectively,
$O(n,p,K) = \big\{ g \in GL(n,K): g^T E_p g = E_p \big\} = \big\{g \in GL(n,K): [gx,gy] = [x,y] \text{ for all } x,y \in X\big\}$ $\big)$,
where $g^T$ is transpose of matrix $g$. Denote by $SO(n,K)$ $\big($respectively, $SO(n,p,K)$ $\big)$ the special orthogonal (respectively, special pseudo-orthogonal) subgroup of $GL(n,K)$, i.e. $SO(n,K) = \big\{g \in O(n,K): \det g = 1\big\}$,
$\big($$SO(n,p,K) = \big\{g \in O(n,p,K): \det g = 1\big\}$$\big)$.
\end{sloppypar}

Now, give the necessary and suficient conditions of $G$-equivalence of strongly regular  $I$-paths $x(t)$ and $y(t)$ by matrices $M(x)(t)$ and $M(y)(t)$, when
 $G$ is one of the following groups: $O(n,K)$, $O(n,p,K)$, $SO(n,K)$, $SO(n,p,K)$.

%%%%%%%%%%%%%%%%%%%%%%%%%%%%%%%%%%%%%%%%%%%%%%%%%%
%%%%%%%%%%%%%%%%%%%%%%%%%%%%%%%%%%%%%%%%%%%%%%%%%%%%%%%%%%%%%%%%%%%%%%%%%%%%%%%%%%%%%%%%%%%%%%%%%%%%
\begin{statement}\label{utv_2_2}
(i). Two strongly regular $I$-paths $x(t)$ and $y(t)$ are $O(n,K)$-equivalent $\big($respectively, $SO(n,K)$-equivalent$\big)$, if and only if the following equalities hold:

\begin{equation}\label{fla_1} \big(M(x)\big)^{-1}(t)M'(x)(t) = \big(M(y)\big)^{-1}M'(y)(t); \end{equation}
\begin{equation}\label{fla_2} M^T(x)(t)M(x)(t) = M^T(y)(t)M(y)(t); \end{equation}

$\big($respectively, the equalities \eqref{fla_1}, \eqref{fla_2} and the equality
\begin{equation}\label{fla_3} \det M(x)(t) = \det M(y)(t) \big) \end{equation}
hold for all $t \in I$.

\begin{sloppypar}
(ii). Two strongly regular $I$-paths $x(t)$ and $y(t)$ are $O(n,p,K)$-equivalent $\big($respectively, $SO(n,p,K)$-equivalent$\big)$, $p \in \{1,..., n-1\}$, if and only if the equalities \eqref{fla_1} and \end{sloppypar}
\begin{equation}\label{fla_4}
M^T(x)(t) E_p M(x)(t) = M^T(y)(t) E_p M(y)(t)
\end{equation}
$\big($respectively, equalities \eqref{fla_1}, \eqref{fla_3}, \eqref{fla_4}$\big)$ hold for every $t \in I$.
\end{statement}
%%%%%%%%%%%%%%%%%%%%%%%%%%%%%%%%%%%%%%%%%%%%%%%%%%%%%%%%%%%%%%%%%%%%%%%%%%%%%%%%%%%%%%%%%%%%%%%%%%%%
%%%%%%%%%%%%%%%%%%%%%%%%%%%%%%%%%%%%%%%%%%%%%%%%%%

The proof of Proposition \ref{utv_2_2} follows the proof of Lemma 3 from \cite{liter_3} $\big($see also Theorem 3 from \cite{liter_6}$\big)$.
For $p =n - 1$,$K = \mathbf{R}$ and for the group $O(n,p,K)$ Proposition \ref{utv_2_2} $(ii)$ is obtained in \cite{liter_5}.

\begin{sloppypar}
By Proposition \ref{utv_2_2}, the equalities for respective rational functions of variables $x_i(t)$, $\mathop{x_i}\limits^{(j)}(t)$, $i,j=\overline{1,n}$ are the nessesery and sufficient condition of $G$-equivalence of strongly regular $I$-paths $x(t)$ and $y(t)$ when $G$ is one of the groups $O(n,K)$, $O(n,p,K)$, $SO(n,K)$, $SO(n,p,K)$. In this respect, it is reasonable to consider the differential field $K \langle x \rangle$ of all rational functions of countable number of variables $x_1, \ldots, x_n, \mathop{x_1}\limits^{(1)}(t), \ldots, \mathop{x_n}\limits^{(1)}(t), \ldots,
\mathop{x_1}\limits^{(r)}(t), \ldots, \mathop{x_n}\limits^{(r)}(t), \ldots$, where differentiation is defined by equality
$d\left(\mathop{x_i}\limits^{(r)}\right) = \mathop{x_i}\limits^{(r+1)}$. Elements from $K \langle x \rangle$ are commonly called $d$-rational functions and the field $K \langle x \rangle$ itself is called $d$-field.\end{sloppypar}

If $\varphi \in K \langle x \rangle$ and $\varphi(gx) = \varphi(x)$ for all $g \in G$, it is said that $d$-rational function $\varphi$ is
$G$-invariant in $X = K^n$ with respect to left action $(g, x) \rightarrow gx$ of subgroup $G \subset GL(n,K)$.

\begin{sloppypar}
Since $\big(M(gx)\big)^{-1}M'(gx) = \big(M(x)\big)^{-1} g^{-1} g M'(x) = \big(M(x)\big)^{-1} M'(x)$, the matrix elements $\big(M(x)\big)^{-1}(t)
\big(M'(x)(t)\big)$ are $G$-invariant $d$-rational functions from $K \langle x \rangle$ for all $t \in I$. Moreover, the elements of matrix
$M^T(x)(t) M(x)(t) = \big(a_{ij}(t)\big)_{i,j=1}^n$ $\big($respectively, $M^T(x)(t) E_p M(x)(t) = \big(b_{ij}(t)\big)_{i,j=1}^n$ $\big)$ have a form
$a_{ij}(t) = \left(\mathop{x}\limits^{(i-1)}(t), \mathop{x}\limits^{(j-1)}(t)\right)$ (respectively, $b_{ij}(t) = \left[\mathop{x}\limits^{(i-1)}(t),
\mathop{x}\limits^{(j-1)}(t)\right]$), $i,j=\overline{1,n}$. If $g \in O(n,K)$ $\big($respectively, $g \in O(n,p,K)$ $\big)$, then
$\left(g\mathop{x}\limits^{(i-1)}(t), g\mathop{x}\limits^{(j-1)}(t)\right) = \left(\mathop{x}\limits^{(i-1)}(t), \mathop{x}\limits^{(j-1)}(t)\right)$
(respectively, $\left[g \mathop{x}\limits^{(i-1)}(t), g \mathop{x}\limits^{(j-1)}(t)\right] = \left[\mathop{x}\limits^{(i-1)}(t),
\mathop{x}\limits^{(j-1)}(t)\right]$) for all $t \in I$.
\end{sloppypar}

\begin{sloppypar}
Thus, the matrix elements from \eqref{fla_1} and \eqref{fla_2} are $O(n,K)$-invariant $d$-rational functions from $K \langle x \rangle$
and the matrix elements from \eqref{fla_1} and \eqref{fla_4} are $O(n,p,K)$-invariant $d$-rational functions from $K \langle x \rangle$. Consequently, for establishment of $O(n,K)$-equivalence $\big($respectiely, $O(n,p,K)$-equivalence$\big)$ of strongly regular $I$-paths $x(t)$ and $y(t)$
it is sufficient to establish equalities $\varphi\big(x(t)\big) = \varphi\big(y(t)\big)$ for all $O(n,K)$-invariant $\big($respectively,
$O(n,p,K)$-invariant$\big)$ $d$-rational functions from $K \langle x \rangle$ and all $t \in I$.
\end{sloppypar}

Let $G$ be a subgroup of $GL(n,K)$ and $K \langle x \rangle^G$ be the set of all $G$-invariant $d$-rational functions from $K \langle x \rangle$. It is known
$\big($\cite{liter_8}, \S 3$\big)$, that $K \langle x \rangle^G$ is a differential subfield of $d$-field $K \langle x \rangle$. A system of elements
$A = \{\varphi_j\}_{j \in J}$ from $K \langle x \rangle^G$ is called system of generators of differential field $K \langle x \rangle^G$, if any element
$\varphi \in K \langle x \rangle^G$ may be constucted from the finite number of elements from the set $A$ using a finite number of operations of the $d$-field $K \langle x \rangle^G$. The following theorem is known $\big($see e.g.\cite{liter_8}, \S 12, Theorems 12.5 and 12.7$\big)$.

%%%%%%%%%%%%%%%%%%%%%%%%%%%%%%%%%%%%%%%%%%%%%%%%%%
%%%%%%%%%%%%%%%%%%%%%%%%%%%%%%%%%%%%%%%%%%%%%%%%%%%%%%%%%%%%%%%%%%%%%%%%%%%%%%%%%%%%%%%%%%%%%%%%%%%%
\begin{theorems}\label{teor_2_3}
\begin{sloppypar}
In $d$-field $K \langle x \rangle^{O(n,K)}$ $\big($respectively, $K \langle x \rangle^{SO(n,K)}$ $\big)$ the following polinomials are generators:
$\left(\mathop{x}\limits^{(j-1)}, \mathop{x}\limits^{(j-1)}\right)$, $j=\overline{1,n}$, $\big($respectively, $\left(\mathop{x}\limits^{(j-1)},
\mathop{x}\limits^{(j-1)}\right)$, $\det M(x)$, $j=\overline{1,(n-1)}$ $\big)$.
\end{sloppypar}
\end{theorems}
%%%%%%%%%%%%%%%%%%%%%%%%%%%%%%%%%%%%%%%%%%%%%%%%%%%%%%%%%%%%%%%%%%%%%%%%%%%%%%%%%%%%%%%%%%%%%%%%%%%%
%%%%%%%%%%%%%%%%%%%%%%%%%%%%%%%%%%%%%%%%%%%%%%%%%%

The following theorem is a variant of Theorem \ref{teor_2_3} for $d$-field $\mathbf{C} \langle x \rangle^{O(n,p,\mathbf{C})}$ and $d$-field
$\mathbf{C} \langle x \rangle^{SO(n,p,\mathbf{C})}$.

%%%%%%%%%%%%%%%%%%%%%%%%%%%%%%%%%%%%%%%%%%%%%%%%%%
%%%%%%%%%%%%%%%%%%%%%%%%%%%%%%%%%%%%%%%%%%%%%%%%%%%%%%%%%%%%%%%%%%%%%%%%%%%%%%%%%%%%%%%%%%%%%%%%%%%%
\begin{theorems}\label{teor_2_4}
\begin{sloppypar}In $d$-field $\mathbf{C} \langle x \rangle^{O(n,p,\mathbf{C})}$ $\big($respectively, $\mathbf{C} \langle x \rangle^{SO(n,p,\mathbf{C})}$ $\big)$ the following polinomials are generators:
$\left[\mathop{x}\limits^{(j-1)}, \mathop{x}\limits^{(j-1)}\right]$, $j=\overline{1,n}$, $\big($respectively, $\left[\mathop{x}\limits^{(j-1)},
\mathop{x}\limits^{(j-1)}\right]$, $i^{(n-p)}\det M(x)$, $j=\overline{1,(n-1)}$, where $i^2 = -1$ $\big)$.\end{sloppypar}
\end{theorems}

\textit{Proof.}
Let $H = \{H_{km}\}_{k,m=1}^n \in GL(n,\mathbf{C})$, where $H_{kk} = 1$ for $k = \overline{1,p}$, $H_{kk} = i$ (here $i^2 = -1$) for $k = \overline{(p+1),n}$,
$H_{km} = 0$ for $k \not= m$. It is clear that $H^T = H$ and $H^2 = E_p = H^{-2}$. Let us show that$G(H) = O(n,\mathbf{C})$, where $G(H) = \big\{HgH^{-1}: g \in O(n,p,\mathbf{C})\big\}$. Indeed, if $g_1 = HgH^{-1}$, where $g \in O(n,p,\mathbf{C})$, then using equalities $H^T H = E_p$, $H^T = E_p H^{-1}$, we have
\begin{multline*}
g_1^T g_1 = (HgH^{-1})^T (HgH^{-1}) = (H^T)^{-1} g^T H^T H g H^{-1} = (H^T)^{-1} g^T E_p g H^{-1} = \\
= (H^T)^{-1} E_p H^{-1} =  (E_p H^{-1})^{-1} E_p H^{-1} = H E_p^{-1} E_p H^{-1} = E,
\end{multline*}
i.e. $g_1 = HgH^{-1} \in O(n,\mathbf{C})$, therefore $G(H) \subset O(n,\mathbf{C})$. The inclusion  $O(n,\mathbf{C}) \subset G(H)$ is established in the same manner.

Now, let us show that $d$-rational function $f(x)$ is invariant with respect to the group $O(n,p,\mathbf{C})$, if and only if $d$-rational function
$\varphi(x) = f(H^{-1}x)$ is invariant with respect to the group $O(n,\mathbf{C})$.

Indeed, if $f \in \mathbf{C} \langle x \rangle^{O(n,p,\mathbf{C})}$, then for $\varphi(x) = f(H^{-1} x)$ and $h \in O(n,\mathbf{C})$, we have
$h = H g H^{-1}$ for some $g \in O(n,p,\mathbf{C})$ and
\[ \varphi(hx) = \varphi(H g H^{-1}x) = f(H^{-1} H g H^{-1} x) = f(g H^{-1} x) = f(H^{-1} x) = \varphi(x), \]
i.e. $\varphi \in \mathbf{C} \langle x \rangle^{O(n,\mathbf{C})}$.

In the same way it is shown that the inclusion $\varphi \in \mathbf{C} \langle x \rangle^{O(n,\mathbf{C})}$ implies the inclusion
$f \in \mathbf{C} \langle x \rangle^{O(n,p,\mathbf{C})}$.

By Theorem \ref{teor_2_3}, polinomials $\varphi_j(x) = \left(\mathop{x}\limits^{(j-1)}, \mathop{x}\limits^{(j-1)}\right)$, $j = \overline{1,n}$,
are the system of generators in  $d$-field $\mathbf{C} \langle x \rangle^{O(n,\mathbf{C})}$. Thus, due to proven above, the system of polinomials
$f_j(x) = \varphi_j(Hx)$, $j = \overline{1,n}$ is the system of generators in $d$-field $\mathbf{C} \langle x \rangle^{O(n,p,\mathbf{C})}$.
Finally, note that
\begin{multline*}
f_j(x) = \left(H\mathop{x}\limits^{(j-1)}, H\mathop{x}\limits^{(j-1)}\right) = \left(H\mathop{x}\limits^{(j-1)}\right)^T
\left(H\mathop{x}\limits^{(j-1)}\right) = \\
= \left(\mathop{x}\limits^{(j-1)}\right)^T H^2 \mathop{x}\limits^{(j-1)} = \left(\mathop{x}\limits^{(j-1)}\right)^T E_p \mathop{x}\limits^{(j-1)} =
\left[\mathop{x}\limits^{(j-1)}, \mathop{x}\limits^{(j-1)}\right], \quad j = \overline{1,n}.
\end{multline*}

Similarly, using Theorem  \ref{teor_2_3} it is established that the polinomials
 $\left[\mathop{x}\limits^{(j-1)}, \mathop{x}\limits^{(j-1)}\right]$, $j = \overline{1,(n-1)}$ together with the polinomial
\[ \det M(Hx) = \det M \big(\{x_1, \ldots, x_p, ix_{p+1}, \ldots, ix_n\}\big) = i^{(n-p)}\det M(x) \] are the system of generators in $d$-field $\mathbf{C} \langle x \rangle^{SO(n,p,\mathbf{C})}$.

\begin{flushright}$\Box$\end{flushright}
%%%%%%%%%%%%%%%%%%%%%%%%%%%%%%%%%%%%%%%%%%%%%%%%%%%%%%%%%%%%%%%%%%%%%%%%%%%%%%%%%%%%%%%%%%%%%%%%%%%%
%%%%%%%%%%%%%%%%%%%%%%%%%%%%%%%%%%%%%%%%%%%%%%%%%%

Combining Proposition \ref{utv_2_2} and Theorems \ref{teor_2_3} and  \ref{teor_2_4} we have the following

%%%%%%%%%%%%%%%%%%%%%%%%%%%%%%%%%%%%%%%%%%%%%%%%%%
%%%%%%%%%%%%%%%%%%%%%%%%%%%%%%%%%%%%%%%%%%%%%%%%%%%%%%%%%%%%%%%%%%%%%%%%%%%%%%%%%%%%%%%%%%%%%%%%%%%%
\begin{consequence}\label{sl_2_5}
\begin{sloppypar}
(i). Two strongly regular $I$-paths $x(t)$ and $y(t)$ are $O(n,K)$-equivalent $\big($respectively, $SO(n,K)$-equivalent$\big)$, if and only if
 $\left(\mathop{x}\limits^{(j-1)}(t), \mathop{x}\limits^{(j-1)}(t)\right) = \left(\mathop{y}\limits^{(j-1)}(t), \mathop{y}\limits^{(j-1)}(t)\right)$
for all $t \in I$, $j = \overline{1,n}$ $\big($respectively,
$\left(\mathop{x}\limits^{(j-1)}(t), \mathop{x}\limits^{(j-1)}(t)\right) = \left(\mathop{y}\limits^{(j-1)}(t), \mathop{y}\limits^{(j-1)}(t)\right)$ and
$\det M(x)(t) = \det M(y)(t)$ for all $t \in I$, $j = \overline{1,(n-1)}$ $\big)$.\end{sloppypar}

\begin{sloppypar}
(ii). Two strongly regular $I$-paths $x(t)$ and $y(t)$ are $O(n,p,K)$-equivalent $\big($respectively, $SO(n,p,K)$-equivalent$\big)$,
$p \in \{1, \ldots, n-1\}$, if and only if
$\left[\mathop{x}\limits^{(j-1)}(t), \mathop{x}\limits^{(j-1)}(t)\right] = \left[\mathop{y}\limits^{(j-1)}(t), \mathop{y}\limits^{(j-1)}(t)\right]$
for all $t \in I$, $j = \overline{1,n}$ $\big($respectively,
$\left[\mathop{x}\limits^{(j-1)}(t), \mathop{x}\limits^{(j-1)}(t)\right] = \left[\mathop{y}\limits^{(j-1)}(t), \mathop{y}\limits^{(j-1)}(t)\right]$ and
\quad $\det M(x)(t) = \det M(y)(t)$ for all $t \in I$, $j = \overline{1,(n-1)}$ $\big)$.
\end{sloppypar}
\end{consequence}
%%%%%%%%%%%%%%%%%%%%%%%%%%%%%%%%%%%%%%%%%%%%%%%%%%%%%%%%%%%%%%%%%%%%%%%%%%%%%%%%%%%%%%%%%%%%%%%%%%%%
%%%%%%%%%%%%%%%%%%%%%%%%%%%%%%%%%%%%%%%%%%%%%%%%%%

%%%%%%%%%%%%%%%%%%%%
%%%%%%%%%%%%%%%%%%%%%%%%%%%%%%%%%%%%%%%%
%%%%%%%%%%%%%%%%%%%%%%%%%%%%%%%%%%%%%%%%%%%%%%%%%%%%%%%%%%%%
%%%%%%%%%%%%%%%%%%%%%%%%%%%%%%%%%%%%%%%%%%%%%%%%%%%%%%%%%%%%%%%%%%%%%%%%%%%%%%%%
\section{Equivalence of paths with respect to the group of movements of pseudo-Euclidean space}

\hspace*{\parindent}Consider the group $Aff(X)$ of all affine transformation of the space $X$. Each transformation from $Aff(X)$ is a superposition of linear nondegenerate transformation $g \in GL(n,K)$ and shift generated by element $u = \{u_i\}_{i=1}^n \in K^n$, i.e. any affine transformation $(u,g) \in Aff(X)$ acts in $X$ by the following rule
\begin{equation}\label{fla_5}
(u,g)(x) = gx + u,
\end{equation}
where $x \in X$, $u \in K^n$, $g \in GL(n,K)$ $\big($here the vector spaces $X$ and $K^n$ are identified$\big)$.

In the group $Aff(X)$ the operation of multiplication is defined by the eqiality
\[ (u,g)(v,h) = (u+gv,gh), \]
where $u,v \in K^n$, $g,h \in GL(n,K)$. In this case, it is said that the group $Aff(X)$ is a semi-direct product of the groups $K^n$ and $GL(n,K)$, that is written as $Aff(X) = K^n \lhd GL(n,K)$. If $G$ is a subgroup of $GL(n,K)$, then the set $K^n \lhd G = \{(u,g) \in K^n \lhd GL(n,K): g \in G\}$ is a subgroup of $K^n \lhd GL(n,K)$, and it is also called semi-direct product of the groups $K^n$ and $G$.

It is known $\big($see e.g. \cite{liter_9}, ch. XVII, \S 2$\big)$, that the group $\mathbf{R}^n \lhd O(n,\mathbf{R})$ coincides with the group of all movements of the Euclidean space
$\big(\mathbf{R}^n, (.\,,.)\big)$, i.e. with the group of all bijections $U$ from $\mathbf{R}^n$ onto $\mathbf{R}^n$, for which $(Ux, Uy) = (x, y)$ for all
$x, y \in \mathbf{R}^n$. Similarly, the group $\mathbf{R}^n \lhd O(n,p,\mathbf{R})$ is the group of all movements of the pseudo-Euclidean space
$\big(\mathbf{R}^n, [.\,,.]\big)$, i.e. the group of all bijections $V$ from
 $\mathbf{R}^n$ onto $\mathbf{R}^n$, for which $[Vx,Vy] = [x,y]$ for all
$x,y \in \mathbf{R}^n$ $\big($see e.g. \cite{liter_10}, ch. III, \S 1$\big)$.

Let $H$ be a subroup of $K^n \lhd GL(n,K)$. Two $I$-paths $x(t)$ and $y(t)$ defined in $X$ are called $H$-equivalent, if there exists $(u,g) \in H$ such that $y(t) = gx(t) + u$ for all $t \in I$.

The following proposition reduces the problem of $K^n \lhd G$-equivalence of $I$-paths $x(t)$ and $y(t)$ to the problem of $G$-equivalence of $I$-paths $x'(t)$ and $y'(t)$.

%%%%%%%%%%%%%%%%%%%%%%%%%%%%%%%%%%%%%%%%%%%%%%%%%%
%%%%%%%%%%%%%%%%%%%%%%%%%%%%%%%%%%%%%%%%%%%%%%%%%%%%%%%%%%%%%%%%%%%%%%%%%%%%%%%%%%%%%%%%%%%%%%%%%%%%
\begin{statement}\label{utv_3_1}
\begin{sloppypar}
Two $I$-paths $x(t)$ and $y(t)$ defined in $X$ are $K^n \lhd G$-equivalent, if and only if $I$-paths $x'(t)$ and $y'(t)$ are $G$-equivalent.
\end{sloppypar}
\end{statement}

\textit{Proof.}
If $I$-paths $x(t)$ and $y(t)$ are $K^n \lhd G$-equivalent, then $y(t) = gx(t) + u$ for all $t \in I$ and some $u = \{u_i\}_{i=1}^n \in K^n$ and
$g = \{g_{ij}\}_{i,j=1}^n \in G$. Since
\[ y_i(t) = \sum_{j=1}^n g_{ij} x_j(t) + u_i, \]
the equation
\begin{equation}\label{fla_8}
y'_i(t) = \sum_{j=1}^n g_{ij} x'_j(t)
\end{equation}
is valid for all $i = \overline{1,n}$, i.e. $y'(t) = g x'(t)$, $t \in I$, that implies the $G$-eqivalence of $I$-paths $x'(t)$ and $y'(t)$.

Conversely, let the equality $y'(t) = g x'(t)$ hold for some $g = \{g_{ij}\}_{i,j=1}^n \in G$ and all $t \in I$. Then the equality \eqref{fla_8} holds, therefore for $u_i(t) = y_i(t) - \sum\limits_{j=1}^n g_{ij} x_j(t)$ we have $u'_i(t) = 0$ for all $t \in I$, i.e. $u_i(t) = u_i^{(0)} \in K$, $t \in I$,
$i = \overline{1,n}$. Consequently, for $u = \{u_i^{(0)}\}_{i=1}^n \in K^n$ the equality $y(t) = gx(t) + u$ holds for all $t \in I$, i.e. $I$-paths $x(t)$ and $y(t)$ are $K^n \lhd G$-equivalent.
\begin{flushright}$\Box$\end{flushright}
%%%%%%%%%%%%%%%%%%%%%%%%%%%%%%%%%%%%%%%%%%%%%%%%%%%%%%%%%%%%%%%%%%%%%%%%%%%%%%%%%%%%%%%%%%%%%%%%%%%%
%%%%%%%%%%%%%%%%%%%%%%%%%%%%%%%%%%%%%%%%%%%%%%%%%%

Proposition \ref{utv_3_1} and Corollary \ref{sl_2_5} imply the following criterion of $K^n \lhd G$-equivalence of $I$-paths $x(t)$ and $y(t)$, in case when
$G$ is one of the groups $O(n,K)$, $O(n,p,K)$, $SO(n,K)$, $SO(n,p,K)$, $p \in \{1, \ldots, n-1\}$.

%%%%%%%%%%%%%%%%%%%%%%%%%%%%%%%%%%%%%%%%%%%%%%%%%%
%%%%%%%%%%%%%%%%%%%%%%%%%%%%%%%%%%%%%%%%%%%%%%%%%%%%%%%%%%%%%%%%%%%%%%%%%%%%%%%%%%%%%%%%%%%%%%%%%%%%
\begin{theorems}\label{teor_3_2}
Let $x(t)$ and $y(t)$ be $I$-paths in $X$, such that $I$-paths $x'(t)$ and $y'(t)$ are strongly regular. Then

\begin{sloppypar}
(i). $I$-paths $x(t)$ and $y(t)$ are $K^n \lhd O(n,K)$-equivalent $\big($respectively, $K^n \lhd O(n,p,K)$-equivalent$\big)$, if and only if
 $\left(\mathop{x}\limits^{(m)}(t), \mathop{x}\limits^{(m)}(t)\right) = \left(\mathop{y}\limits^{(m)}(t), \mathop{y}\limits^{(m)}(t)\right)$
$\big($respectively, $\left[\mathop{x}\limits^{(m)}(t), \mathop{x}\limits^{(m)}(t)\right] =
\left[\mathop{y}\limits^{(m)}(t), \mathop{y}\limits^{(m)}(t)\right]$ $\big)$ for all $t \in I$, $m = \overline{1,n}$.
\end{sloppypar}

\begin{sloppypar}
(ii). $I$-paths $x(t)$ and $y(t)$ are $K^n \lhd SO(n,K)$-equivalent $\big($respectively, $K^n \lhd SO(n,p,K)$-equivalent$\big)$, if and only if $\left(\mathop{x}\limits^{(m)}(t), \mathop{x}\limits^{(m)}(t)\right) = \left(\mathop{y}\limits^{(m)}(t), \mathop{y}\limits^{(m)}(t)\right)$
and $\det M'(x)(t) = \det M'(y)(t)$ $\big($respectively, $\left[\mathop{x}\limits^{(m)}(t), \mathop{x}\limits^{(m)}(t)\right] =
\left[\mathop{y}\limits^{(m)}(t), \mathop{y}\limits^{(m)}(t)\right]$ and $\det M'(x)(t) = \det M'(y)(t)$ $\big)$ for all $t \in I$, $m = \overline{1,(n-1)}$.
\end{sloppypar}
\end{theorems}
%%%%%%%%%%%%%%%%%%%%%%%%%%%%%%%%%%%%%%%%%%%%%%%%%%%%%%%%%%%%%%%%%%%%%%%%%%%%%%%%%%%%%%%%%%%%%%%%%%%%
%%%%%%%%%%%%%%%%%%%%%%%%%%%%%%%%%%%%%%%%%%%%%%%%%%

%%%%%%%%%%%%%%%%%%%%
%%%%%%%%%%%%%%%%%%%%%%%%%%%%%%%%%%%%%%%%
%%%%%%%%%%%%%%%%%%%%%%%%%%%%%%%%%%%%%%%%%%%%%%%%%%%%%%%%%%%%
%%%%%%%%%%%%%%%%%%%%%%%%%%%%%%%%%%%%%%%%%%%%%%%%%%%%%%%%%%%%%%%%%%%%%%%%%%%%%%%%
\section{Equivalence of curves with respect to the action of the groups $K^n \lhd O(n,p,K)$, $K^n \lhd SO(n,p,K)$, $O(n,p,K)$ and $SO(n,p,K)$.}

\hspace*{\parindent}Two paths $x: I_1 \rightarrow X$ and $y: I_2 \rightarrow X$ are called $D$-equivalent (respectively, $D_+$-equivalent)
$\big($\cite{liter_8}, \S 4$\big)$, if there exists $C^\infty$-diffeomorphism $\varphi$ from $I_2$ onto $I_1$ such that $\varphi'(t) \not= 0$
$\big($respectively, $\varphi'(t) > 0$ $\big)$ and $y(t) = x\big(\varphi(t)\big)$ for all $t \in I_2$.

Obviously, $D$-equivalence (respectively, $D_+$-equivalence) is a relation of equivalence on the set of all paths defined in $X$.

A class $\gamma = \overline{x}$ (respectively, $\gamma = \hat{x}$) of all paths in $X$,  $D$-eqiuvalent (respectively, $D_+$-equivalent) for the path $x(t)$, is called curve (respectively, oriented curve) in $X$ generated by this path $\big($\cite{liter_8}, \S 4$\big)$. In this case, a path $y$ from the class $\gamma$ is called parametrization of the curve $\gamma$.

For every $I$-path $x: I \rightarrow X$ its image $x(I)$ in $X$ is called support of path $x$ and denoted by $\widetilde{x}$, i.e. $\widetilde{x} =
\big\{\{x_j(t)\}_{j=1}^n: t \in I\big\} \subset X$.

It is clear that $D$-equivalent paths $x(t)$ and $y(t)$ have the same support, i.e. $\widetilde{x} = \widetilde{y}$. Therefore, the set $\widetilde{x}$ is said to be the support of the curve
$\gamma = \overline{x}$.

Note, that, in general, the equality of supports $\widetilde{x} =
\widetilde{y}$ of paths $x(t)$ and $y(t)$, does not imply their
$D$-equivalence. At the same time, according to Theorem
\ref{teor_2_1}, for simply-regular paths $x(t)$ and $y(t)$ the equality
$\widetilde{x} = \widetilde{y}$ implies the $D$-equivalence of paths
$x(t)$ and $y(t)$.

Thus, up to $ D $-equivalence, simple regular paths are uniquely determined by their supports, and therefore, the curves generated by such paths may be identified with their supports.

If $x$ is a simple (respectively, regular) path and $y \in
\overline{x}$, then, obviously, the path $y$ is also simple (respectively,
regular). The curve $\gamma$ is called simple (respectively,
regular, strongly-regular), if $y = \overline{x}$, where $x(t)$
is a simple (respectively, regular, strongly-regular) path.

Let $H$ be a subgroup of the group $Aff(X)$. If $h = (u,g) \in
H$ and $x(t)$ is a $I$-path in $X$, then $y(t) = h\big(x(t)\big)
= g\big(x(t)\big) + u$ $\big($see equality \eqref{fla_5}$\big)$ is a
$I$-path in $X$. The curve (oriented curve) generated by $I$-path
$y(t)$ will be denoted by $h\gamma$, where $\gamma =
\overline{x}$ (respectively, $\gamma = \hat{x}$), i.e. $h\gamma =
\overline{hx}$ (respectively, $h\gamma = \widehat{hx}$).

%%%%%%%%%%%%%%%%%%%%%%%%%%%%%%%%%%%%%%%%%%%%%%%%%%
%%%%%%%%%%%%%%%%%%%%%%%%%%%%%%%%%%%%%%%%%%%%%%%%%%%%%%%%%%%%%%%%%%%%%%%%%%%%%%%%%%%%%%%%%%%%%%%%%%%%
\begin{statement}\label{utv_4_1}
Let $\gamma_1 = \overline{x}$, $\gamma_2 = \overline{y}$ be simple
regular curves generated by $I_1$-path $x(t)$ and $I_2$-path
$y(t)$ respectively. The following conditions are equivalent:

(i) there exists $h \in H$ such that $\gamma_2 = h\gamma_1$;

(ii) there exist $h \in H$ and $C^\infty$-diffeomorphism
$\varphi: I_2 \rightarrow I_1$ such that $\varphi'(t) \not= 0$ and $y(t)
= h x\big(\varphi(t)\big)$ for all $t \in I_2$.
\end{statement}

\textit{Proof.}$(i) \Rightarrow (ii)$. Equality $\overline{y} =
\gamma_2 = h\overline{x}$ means that $\widetilde{y} =
h\widetilde{x}$, where $h = (u,g)$, $u = \{u_i\}_{i=1}^n \in K^n$,
$g = \{g_{ij}\}_{i,j=1}^n \in GL(n,K)$, and $z(t) = (hx)(t) = gx(t)
+ u$ $\big($see equality \eqref{fla_5}$\big)$ is a simple $I_1$-path
from the class $\overline{y}$.

Hence, supports of simply-regular paths $y(t) = \big\{y_1(t),
\ldots, y_n(t)\big\}$ and \\$z(t) = \left\{ \sum\limits_{j=1}^n g_{1j}
x_j(t) + u_1, \ldots, \sum\limits_{j=1}^n g_{nj} x_j(t) +
u_n\right\}$ coincide. By Theorem
\ref{teor_2_1}, there exists  a $C^\infty$-diffeomorphism
$\varphi: I_2 \rightarrow I_1$ such that $\varphi'(t) \not= 0$ and $y(t)
= z\big(\varphi(t)\big) = h x\big(\varphi(t)\big)$ for all $t \in
I_2$.

$(ii) \Rightarrow (i)$. Let $y(t) = hx\big(\varphi(t)\big)$ for all
$t \in I_2$ and some $h \in H$ and $C^\infty$-diffeomorphism
$\varphi: I_2 \rightarrow I_1$ with $\varphi'(t) \not= 0$, $t \in I_2$.
Paths $x(t)$ and $x\big(\varphi(t)\big)$ are $D$-equivalent, and
thus $\overline{x} = \overline{x\circ\varphi}$. Consequently,
$\gamma_2 = \overline{y} = \overline{h(x\circ\varphi)} =
h\left(\overline{x\circ\varphi}\right) = h\overline{x} = h\gamma_1$.
\begin{flushright}$\Box$\end{flushright}
%%%%%%%%%%%%%%%%%%%%%%%%%%%%%%%%%%%%%%%%%%%%%%%%%%%%%%%%%%%%%%%%%%%%%%%%%%%%%%%%%%%%%%%%%%%%%%%%%%%%
%%%%%%%%%%%%%%%%%%%%%%%%%%%%%%%%%%%%%%%%%%%%%%%%%%

Let $H$ be a subgroup of the group $Aff(X) = K^n \lhd GL(n,K)$. Curves
$\gamma_1$ and $\gamma_2$ are called $H$-equivalent, if there
exists $h \in H$ such that $\gamma_2 = h \gamma_1$.

Clearly, the $H$-equivalence of paths $x(t)$ and $y(t)$ implies the
$H$-equivalence of curves $\overline{x}$ and $\overline{y}$. The
converse statement, in general, does not hold (see Proposition
\ref{utv_4_1}).

Of course, there is the problem of existence of parameterizations $x(t)$ and $y(t)$ of curves $\gamma_1$ and $\gamma_2$, respectively, such that $H$-equivalence of curves $\gamma_1$
and $\gamma_2$ implies $H$-equivalence of paths $x(t)$
and $y(t)$. In the case of oriented curves for the groups
$\mathbf{R}^n \lhd O(n,\mathbf{R})$ and $\mathbf{R}^n \lhd
SO(n,\mathbf{R})$ this problem is solved in \cite{liter_4} using invariant parameterizations,
and for the group $\mathbf{R}^n \lhd SL(n,\mathbf{R})$ in \cite{liter_2}, and for
the centrally-affine groups  in \cite{liter_1}. In
this section there is given the solution to the above mentioned problem for
the oriented curves wth respect to the action of the groups $K^n \lhd
O(n,p,K)$, $K^n \lhd SO(n,p,K)$, $O(n,p,K)$ and $SO(n,p,K)$.

Let $x(t)$ be an $I$-path in $X = K^n$ such that the inequation $\big[x'(t),
x'(t)\big] \not= 0$ holds for all $t \in I = (a,b)$ $\big($such
paths will be called non-degenerate and, respectively, the curve
$\overline{x}$ (oriented curve $\hat{x}$) will also be called
non-degenerate$\big)$. In this case, it is clear that $x'(t) \not= 0$
for all $t \in I$, i.e. non-degenerate $I$-path $x(t)$ is always
regular. Note, that in \cite{liter_4} the term <<non-degenerate
path>> means strongly regular path in our sense.

Since $\big[x'(t), x'(t)\big]$ is a continuous function on $I$,
then for any $c,d \in (a,b)$, $c < d$, there exists a finite
integral
\begin{equation}\label{fla_9}
l_x(c,d) = \int_c^d \Big|\big[x'(t),x'(t)\big]\big|^{1/2}dt.
\end{equation}
Therefore, there exist finite or infinite limits
\begin{equation}\label{fla_10}
l_x(a,d) = \lim_{c \to a}l_x(c,d),
\end{equation}
\begin{equation}\label{fla_11}
l_x(c,b) = \lim_{d \to b}l_x(c,b).
\end{equation}

As in \cite{liter_2}, \cite{liter_4} we say that $I$-path
$x(t)$ has a type

(L1), if $l_x(a,d) < \infty$ and $l_x(c,b) < \infty$;

(L2), if $l_x(a,d) < \infty$ and $l_x(c,b) = +\infty$;

(L3), if $l_x(a,d) = +\infty$ and $l_x(c,b) < \infty$;

(L4), if $l_x(a,d) = +\infty$ and $l_x(c,b) = +\infty$.

\noindent Clearly, the definition of the type of a non-degenerate path
$x(t)$ is independent of choice of points $c,d \in I$.

Define interval $I(x) = \big(A(x), B(x)\big) \subset \mathbf{R}$ by the following rule:
\begin{itemize}
\item[(i)] if $x(t)$ has type (L1), then $A(x) = 0$ and $B(x) = \int\limits_a^b \big|\big[x'(t),x'(t)\big]\big|^{1/2} dt$;
\item[(ii)] if $x(t)$ has type (L2), then $A(x) = 0$ and $B(x) = +\infty$;
\item[(iii)] if $x(t)$ has type (L3), then $A(x) = -\infty$ and $B(x) = 0$;
\item[(iv)] if $x(t)$ has type (L4), then $A(x) = -\infty$ and $B(x) = +\infty$;
\end{itemize}

Let us specify the connection between the types of  $D_+$-equivalent
non-degenerate paths. Let $x(t)$ be an
$I_1$-path and let $y(t)$ be an $I_2$-path, $I_1 = (a_1,b_1)$, $I_2 =
(a_2,b_2)$. Assume, that paths $x(t)$ and $y(t)$ are non-degenerate
and $D_+$-equivalent, in particular, there exists a
$C^\infty$-diffeomorphism from $I_2$ onto $I_1$ such that
$\varphi'(t)
> 0$ and $y(t) = x\big(\varphi(t)\big)$ for all $t \in I_2$. Note, if the paths $x(t)$ and $y(t)$ are $D_+$-equivalent the
non-degeneracy of the path $x(t)$ implies non-degeneracy of the path
$y(t)$. Due to \eqref{fla_9}, for $a_2 < c < d < b_2$ we have
 \[ l_y(c,d) = \int_c^d \big|\big[y'(t),y'(t)\big]\big|^{1/2}dt
= \int_c^d
\left|\frac{d\varphi}{dt}\right|\cdot\left|\left[\frac{d}{d\varphi}x\big(\varphi(t)\big),\frac{d}{d\varphi}x\big(\varphi(t)\big)\right]\right|^{1/2}dt
\]

Since $\varphi'(t) > 0$, it follows that
\[ l_y(c,d) = \int_{\varphi(c)}^{\varphi(d)} \left|\left[\frac{d}{d\varphi}x\big(\varphi(t)\big),\frac{d}{d\varphi}x\big(\varphi(t)\big)\right]\right|^{1/2}d\varphi = l_x\big(\varphi(c),\varphi(d)\big). \]

By equalities \eqref{fla_10} and
$\eqref{fla_11}$ we have
\begin{equation}\label{fla_12}
l_y(a_2,d) = \lim_{c \to a_2}l_y(c,d) = \lim_{s \to a_1}l_x\big(s,\varphi(d)\big) = l_x\big(a_1,\varphi(d)\big),
\end{equation}
\begin{equation}\label{fla_13}
l_y(c,b_2) = \lim_{d \to b_2}l_y(c,d) = \lim_{r \to b_1}l_x\big(\varphi(c),r\big) = l_x\big(\varphi(c), b_1\big),
\end{equation}

It means that $D_+$-equivalent paths $x(t)$ and $y(t)$
have same type.

Equalities \eqref{fla_12} and \eqref{fla_13} imply the following

%%%%%%%%%%%%%%%%%%%%%%%%%%%%%%%%%%%%%%%%%%%%%%%%%%
%%%%%%%%%%%%%%%%%%%%%%%%%%%%%%%%%%%%%%%%%%%%%%%%%%%%%%%%%%%%%%%%%%%%%%%%%%%%%%%%%%%%%%%%%%%%%%%%%%%%
\begin{statement}\label{utv_4_2}
If non-degenerate $I_1$-path $x(t)$ and $I_2$-path $y(t)$ are
$D_+$-equivalent, then the paths $x(t)$ and $y(t)$ have same
type and $I(x) = I(y)$.
\end{statement}
%%%%%%%%%%%%%%%%%%%%%%%%%%%%%%%%%%%%%%%%%%%%%%%%%%%%%%%%%%%%%%%%%%%%%%%%%%%%%%%%%%%%%%%%%%%%%%%%%%%%
%%%%%%%%%%%%%%%%%%%%%%%%%%%%%%%%%%%%%%%%%%%%%%%%%%

Also note the following useful property of $K^n \lhd O(n,p,K)$-equivalent
non-degenerate paths.

%%%%%%%%%%%%%%%%%%%%%%%%%%%%%%%%%%%%%%%%%%%%%%%%%%
%%%%%%%%%%%%%%%%%%%%%%%%%%%%%%%%%%%%%%%%%%%%%%%%%%%%%%%%%%%%%%%%%%%%%%%%%%%%%%%%%%%%%%%%%%%%%%%%%%%%
\begin{statement}\label{utv_4_3}
Let $x(t)$ be a non-degenerate $I$-path in $K^n$, $(u,g) \in K^n
\lhd O(n,p,K)$ and $y(t) = gx(t) + u$, $t \in I$. Then the paths
$x(t)$ and $y(t)$ have same type and $I(x) = I(y)$.
\end{statement}
%%%%%%%%%%%%%%%%%%%%%%%%%%%%%%%%%%%%%%%%%%%%%%%%%%%%%%%%%%%%%%%%%%%%%%%%%%%%%%%%%%%%%%%%%%%%%%%%%%%%
%%%%%%%%%%%%%%%%%%%%%%%%%%%%%%%%%%%%%%%%%%%%%%%%%%

The proof directly follows from the equalities
\[ \big[y'(t),y'(t)\big] = \big[gx'(t),gx'(t)\big] = \big[x'(t),x'(t)\big]. \]

Now define a special parametrization for a non-degenerate
$I$-path $x(t)$. Consider the function $p_x$ from $I=(a,b)$ onto $I(x)$,
defined by the following rule $\big($see \cite{liter_2},
\cite{liter_4}$\big)$: if $I$-path $x(t)$ has

\begin{itemize}
\item[(i)] type (L1) or (L2), then
\[ p_x(t) = l_x(a,t) = \int_a^t \big|\big[x'(t),x'(t)\big]\big|^{1/2}dt; \]
\item[(ii)] type (L3), then
\[ p_x(t) = - l_x(t,b) = - \int_t^b \big|\big[x'(t),x'(t)\big]\big|^{1/2}dt; \]
\item[(iii)] type (L4), then, fixing a point $a_I \in I$, we set
\[ p_x(t) = p_{x,a_I} = l_x(a_I,t) = \int_{a_I}^t \big|\big[x'(t),x'(t)\big]\big|^{1/2}dt. \]
\end{itemize}
Obviously, if we fix another point $a'_I \in I$, we
have that $p_{x,a_I}(t) - p_{x,a'_I}(t) = l_x(a_I,a'_I) = const$,
i.e. for the type (L4), the function $p_x(t)$ is defined up to a
constant. If $I = (-\infty, +\infty)$, then we will always set $a_I
= 0$.

Since $x(t)$ is the non-degenerate $I$-path, then
$\big|[x'(t),x'(t)]\big|
> 0$ for all $t \in I$, and therefore the function $p_x(t)$ is a
$C^\infty$-diffeomorphism.

Due to inequality $p'_x(t) > 0$, $t \in I$, for the function $p_x(t)$
there exists an inverse function $q_x(s)$ from $I(x)$ onto $I$, in addition, $q_x$ is also a $C^\infty$-diffeomorphism and $q'_x(s)
> 0$ for all $s \in I(x)$. Consequently, $y(s) =
x\big(q_x(s)\big)$ is a non-degenerate $I(x)$-path
$D_+$-equivalent to the $I$-path $x(t)$, in particular, $\hat{x} = \hat{y}$.

The following proposition describes useful properties of
parametrization $x\big(q_x(s)\big)$ for the oriented curve
$\hat{x}$.

%%%%%%%%%%%%%%%%%%%%%%%%%%%%%%%%%%%%%%%%%%%%%%%%%%
%%%%%%%%%%%%%%%%%%%%%%%%%%%%%%%%%%%%%%%%%%%%%%%%%%%%%%%%%%%%%%%%%%%%%%%%%%%%%%%%%%%%%%%%%%%%%%%%%%%%
\begin{statement}\label{utv_4_4}
Let $x(t)$ be a non-degenerate $I$-path in $K^n$, $I = (a,b)$,
$\varphi$ be a $C^\infty$-diffeomorphism from $J = (c,d)$ onto $I$,
$\varphi'(r)
> 0$, $r \in J$, $h=(a,g) \in K^n \lhd O(n,p,K)$, $p \in \{1,
\ldots, n-1\}$. Then

\begin{sloppypar}
(i) $p_{hx}(t) = p_x(t)$ $\big($respectively, $q_{hx}(s) =
q_x(s)$ $\big)$ for all $t \in I$ $\big($respectively, for all $s \in
I(x) = I(hx)$ \emph{(}see Proposition
\ref{utv_4_3}\emph{)}$\big)$;\end{sloppypar}

(ii) If $I$-path $x(t)$ has type (L1), (L2) or (L3)
$\big($respectively, \emph{(}L4\emph{)}$\big)$, then
\begin{equation}\label{fla_14}
p_{x\circ\varphi}(r) = p_x\big(\varphi(r)\big)
\end{equation}
and
\begin{equation}\label{fla_15}
\varphi\big(q_{x\circ\varphi}(s)\big) = q_x(s)
\end{equation}
(respectively,
\begin{equation}\label{fla_16}
p_{x\circ\varphi}(r) = p_x\big(\varphi(r)\big) + s_0
\end{equation}
and
\begin{equation}\label{fla_17}
\varphi\big(q_{x\circ\varphi}(s + s_0)\big) = q_x(s)\big)
\end{equation}
for all $r \in J$, $s \in I(x) = I(x\circ\varphi)$ (see Proposition
\ref{utv_4_2}), where $s_0 = l_x\big(\varphi(a_J),a_I\big)$;

(iii) for $I(x)$-path $z(s) = x\big(q_x(s)\big)$ the
equalities $I(x) = I(z)$ and
\begin{equation}\label{fla_18}
\big|[z'(s),z'(s)]\big| = 1, \quad p_z(s) = s = q_z(s)
\end{equation}
hold for all $s \in I(x)$.
\end{statement}

\textit{Proof.} $(i)$ follows from the equalities
\[ \big|[(hx)'(t),(hx)'(t)]\big| = \big|[gx'(t),gx'(t)]\big| = \big|[x'(t),x'(t)]\big|. \]

$(ii)$ Let $I$-path $x(t)$ has type (L1) or (L2). Using
inequality $\varphi'(r)
> 0$, $r \in J$, and, by changing variables $t = \varphi(r)$,
we have that
\begin{multline*}
p_{x\circ\varphi}(r) = \int_c^r \big|[(x\circ\varphi)'(s),(x\circ\varphi)'(s)]\big|^{1/2}ds = \\
=\int_c^r \left|\left[\frac{d\varphi}{ds}\frac{d}{d\varphi}(x\circ\varphi)(s),\frac{d\varphi}{ds}\frac{d}{d\varphi}(x\circ\varphi)(s)\right]\right|^{1/2}ds=\\
= \int_c^r \frac{d\varphi}{ds}\left|\left[\frac{d}{d\varphi}x\big(\varphi(s)\big),\frac{d}{d\varphi}x\big(\varphi(s)\big)\right]\right|^{1/2}ds = \\
= \int_{\varphi(c)}^{\varphi(r)} \big|[x'(t),x'(t)]\big|^{1/2}dt = p_x\big(\varphi(r)\big)
\end{multline*}
for all $r \in J$. Similarly, we obtain the equality
\eqref{fla_14} in case when $x(t)$ is an $I$-path has type (L3).

If $I$-path $x(t)$ has type (L4), then
\[p_{x\circ\varphi}(r) = \int_{a_J}^r
\big|[(x\circ\varphi)'(s),(x\circ\varphi)'(s)]\big|^{1/2}ds, \]
where $a_J \in J$. Therefore,

\begin{multline*}
p_{x\circ\varphi}(r) = \\
= \int_{a_J}^r \frac{d\varphi}{ds}\left|\left[\frac{d}{d\varphi}x\big(\varphi(s)\big),\frac{d}{d\varphi}x\big(\varphi(s)\big)\right]\right|^{1/2}ds
= \int_{\varphi(a_J)}^{\varphi(r)} \big|[x'(t),x'(t)]\big|^{1/2}dt = \\
= \int_{a_I}^{\varphi(r)} \big|[x'(t),x'(t)]\big|^{1/2}dt + \int_{\varphi(a_J)}^{a_I} \big|[x'(t),x'(t)]\big|^{1/2}dt = \\
= p_x\big(\varphi(r)\big) + s_0.
\end{multline*}

%\textbabygamma

Since $q_x: I(x) \rightarrow I$ and $q_{x\circ\varphi}: I(x)
\rightarrow J$ are inverse to the functions $p_x$ and
$p_{x\circ\varphi}$, respectively, the equalities
\eqref{fla_15} and \eqref{fla_17} follow from equalities \eqref{fla_14}
and \eqref{fla_16}, respectively.

$(iii)$ If $I$-path $x(t)$ has type (L1) or (L2), then,
according to $(ii)$, for $\varphi(s) = q_x(s)$, $s \in I(x)$, we
have that
\[ \int_{A(x)}^s \big|[z'(r),z'(r)]\big|^{1/2}dr = p_z(s) = p_{x\circ\varphi}(s) = p_x\big(q_x(s)\big) = s, \]
and hence $\big|[z'(s),z'(s)]\big| = 1$ for all $s \in I(x)$. If
$I$-path $x(t)$ has type (L3), then
\[ - \int_s^{B(x)} \big|[z'(r),z'(r)]\big|^{1/2}dr = p_z(s) = p_x\big(q_x(s)\big) = s, \]
that also implies the equality $\big|[z'(s),z'(s)]\big| = 1$.

If $I$-path $x(t)$ has type (L4), then $I(x) = (-\infty,
+\infty)$, $a_{I(x)} = 0$, and for $\varphi(s) = q_x(s)$, $s \in I(x)
= (-\infty, +\infty)$ we have
\[ \int_0^s \big|[z'(r),z'(r)]\big|^{1/2}dr = p_z(s) = p_{x\circ\varphi}(s) = p_x\big(q_x(s)\big) + l_x\big(q_x(0), a_I\big) = s + s_0 \]
for all $s \in I(x)$, where $s_0 = l_x\big(q_x(0), a_I\big)$ is a
constant. For $s=0$ we have that $s_0 = 0$, i.e. $p_z(s) = s$ for
all $s \in I(x)$. Thus, equality $\big|[z'(s),z'(s)]\big| = 1$ holds for all $s \in I(x)$ in this case too.

Since $q_z(s)$ is an inverse function for $p_z(s) = s$, we have
$q_z(s) = s$. Equality $I(z) = I(x)$ immediately follows from
\eqref{fla_18}.
\begin{flushright}$\Box$\end{flushright}
%%%%%%%%%%%%%%%%%%%%%%%%%%%%%%%%%%%%%%%%%%%%%%%%%%%%%%%%%%%%%%%%%%%%%%%%%%%%%%%%%%%%%%%%%%%%%%%%%%%%
%%%%%%%%%%%%%%%%%%%%%%%%%%%%%%%%%%%%%%%%%%%%%%%%%%

Let $x(t)$ and $y(t)$ be non-degenerate $I$-paths, $h \in K^n \lhd
O(n,p,K)$ and $y(t) = hx(t)$, $t \in I$. From Propositions
\ref{utv_4_3} and \ref{utv_4_4} $(i)$ it follows, that $K^n \lhd
O(n,p,K)$-equivalent non-degenerate curves $\overline{x}$ and
$\overline{y}$ have same parametrization $q_x(s) = q_y(s)$, $s
\in I(x) = I(y)$, i.e. $x \circ q_x \in \overline{x}$ and $y \circ
q_x \in \overline{y}$. In this respect, parametrization of the form
$z(s) = x\big(q_x(s)\big)$ for the curve $\overline{x}$ is obviously to be called invariant parametrization with respect to the action of group
$K^n \lhd O(n,p,K)$.

Propositions \ref{utv_4_2} and \ref{utv_4_4} $(ii)$ also imply
the following

%%%%%%%%%%%%%%%%%%%%%%%%%%%%%%%%%%%%%%%%%%%%%%%%%%
%%%%%%%%%%%%%%%%%%%%%%%%%%%%%%%%%%%%%%%%%%%%%%%%%%%%%%%%%%%%%%%%%%%%%%%%%%%%%%%%%%%%%%%%%%%%%%%%%%%%
\begin{consequence}\label{sl_4_5}
If $x$ is a non-degenerate $I$-path, $y$ is a non-degenerate
$J$-path and $\hat{x} = \hat{y}$, then paths $x$ and $y$ have
same type, $I(x) = I(y)$ and $x\big(q_x(s)\big) = y\big(q_y(s +
s_0)\big)$ for all $s \in I(x)$, where $s_0$ is some constant,
which is equal to zero, in cases, when $I$-path $x(t)$ has one of the types
(L1), (L2) or (L3).
\end{consequence}

\textit{Proof.}
Oriented curve $\hat{x} = \hat{y}$ is the class of all
paths which is $D_+$-equivalent to the path $x(t)$. Therefore,
paths $x$ and $y$ are $D_+$-equivalent, i.e. there exists a
$C^\infty$-diffeomorphism $\varphi: J \rightarrow I$ such that
$\varphi'(r) > 0$ and $y(r) = x\big(\varphi(r)\big)$ for all $r \in
J$. By Proposition \ref{utv_4_2}, the types of paths $x$ and $y$ are same and $I(x) = I(y)$. If
$I$-path $x(t)$ has type (L1), (L2) or (L3), then, by Proposition \ref{utv_4_4} $(ii)$, we have
\[ y\big(q_y(s)\big) = (x\circ\varphi)\big(q_y(s)\big) = x\big(\varphi\big(q_{x\circ\varphi}(s)\big)\big) = x\big(q_x(s)\big) \]
for all $s \in I(x) = I(y)$.

If $I$-path $x(t)$ has type (L4), then again using Proposition \ref{utv_4_4} $(ii)$, we have
\[ y\big(q_y(s)\big) = x\big(\varphi\big(q_{x\circ\varphi}(s)\big)\big) = x\big(q_x(s - s_0)\big) \]
for all $s \in I(x)$, where $s_0 = l_x\big(\varphi(a_J),a_I\big)$.
\begin{flushright}$\Box$\end{flushright}
%%%%%%%%%%%%%%%%%%%%%%%%%%%%%%%%%%%%%%%%%%%%%%%%%%%%%%%%%%%%%%%%%%%%%%%%%%%%%%%%%%%%%%%%%%%%%%%%%%%%
%%%%%%%%%%%%%%%%%%%%%%%%%%%%%%%%%%%%%%%%%%%%%%%%%%

According to Corollary \ref{sl_4_5}, we say that a
non-degenerate oriented curve $\gamma$ has type (L1)
$\big($respectively, (L2), (L3), (L4)$\big)$, if the path
$x\in\gamma$ has the same type; in this case we set $I(\gamma) =
I(x)$.

From Proposition \ref{utv_4_3} it follows, that $K^n \lhd
O(n,p,K)$-equivalent oriented curves $\gamma$ and $\beta$ always
have same types and $I(\gamma) = I(\beta)$.

%%%%%%%%%%%%%%%%%%%%%%%%%%%%%%%%%%%%%%%%%%%%%%%%%%
%%%%%%%%%%%%%%%%%%%%%%%%%%%%%%%%%%%%%%%%%%%%%%%%%%%%%%%%%%%%%%%%%%%%%%%%%%%%%%%%%%%%%%%%%%%%%%%%%%%%
\begin{theorems}\label{teor_4_6}
Let $x$ be a non-degenerate $I$-path, $y$ be a non-degenerate
$J$-path, $x_1(s) = x\big(q_x(s)\big)$, $s \in I(x)$, $y_1(r) =
y\big(q_y(r)\big)$, $r \in I(y)$, $\gamma = \hat{x}_1$, $\beta =
\hat{y}_1$, $h = (u,g) \in K^n \lhd O(n,p,K)$. If the path $x(t)$ has
type (L1), (L2) or $(L3)$, then $\beta = h\gamma$, if and only if
$I(x) = I(y)$ and $y_1(s) = hx_1(s)$ for all $s \in I(x)$.

If the path $x(t)$ has type (L4), then $\beta = h\gamma$,
if and only if $I(x) = I(y) = (-\infty, +\infty)$ and $y_1(s+s_0) =
(h x_1)(s)$ for all $s \in I(x)$ and some $s_0 \in I(x)$.
\end{theorems}

\textit{Proof.}
If the path $x(t)$ has type (L1), (L2) or
(L3), $I(x) = I(y)$ and $y_1(s) = h x_1(s)$ for all $s \in I(x)$,
then, obviously, $\beta = h\gamma$. If the path
$x(t)$ has type (L4) and $y_1(s + s_0) = h x_1(s)$ for all $s
\in I(x)$ and some $s_0 \in I(x)$, then the path $z(s) = y_1(s+s_0)$
also has type (L4) and the inclusion $z\in\beta$ implies the equality $\beta = h\gamma$.

Now, let $\hat{y}_1 = \beta = h\gamma = h\hat{x}_1$. Thus, the
$I(x)$-path $z(s) = (hx_1)(s) = gx_1(s) + u$ is a path from the
equivalence class $\hat{y}_1$. It means that $I(x)$-path $z$ and
$I(y)$-path $y_1$ are $D_+$-equivalent, i.e. there exists a
$C^\infty$-diffeomorphism $\varphi$ from $I(y)$ onto $I(x)$ such that
$y_1(s) = z\big(\varphi(s)\big)$, $\varphi'(s) > 0$ for all $s \in
I(y)$.

By Proposition \ref{utv_4_4} $(iii)$ and Corollary
$\ref{sl_4_5}$, we have that types of the paths $z$ and $y_1$ are
same, $I(z) = I(y_1)$ and
\begin{equation}\label{fla_19}
y_1(s+s_0) = y_1\big(q_{y_1}(s+s_0)\big) = z\big(q_z(s)\big)
\end{equation}
for all $s \in I(z)$, where $s_0$ is a constant, which is equal to zero,
when $I(x)$-path $z(s)$ has type (L1), (L2) or (L3).

According to Proposition \ref{utv_4_4} $(i)$, $(iii)$, we have
 $I(y_1) = I(z) = I(x)$ and $q_z(s) = q_{x_1}(s)$, in
particular,
\begin{equation}\label{fla_20}
z\big(q_z(s)\big) = \big(hx_1\big)\big(q_{x_1}(s)\big) = (hx_1)(s)
\end{equation}
for all $s \in I(x)$.

Equalities \eqref{fla_19} and \eqref{fla_20} imply that
$I(y) = I(y_1) = I(x_1) = I(x)$ and
\begin{equation}\label{fla_21}
y_1(s+s_0) = (hx_1)(s)
\end{equation}
for all $s \in I(x_1)$. If $I$-path $x$ has type (L1), (L2)
or(L3), then $s_0 = 0$, i.e. $y_1(s) = (hx_1)(s)$. If $I$-path $x$
has type (L4), then $I(x) = I(y) = (-\infty, +\infty)$ and equality \eqref{fla_21} holds for some $s_0\in\mathbf{R}$.
\begin{flushright}$\Box$\end{flushright}
%%%%%%%%%%%%%%%%%%%%%%%%%%%%%%%%%%%%%%%%%%%%%%%%%%%%%%%%%%%%%%%%%%%%%%%%%%%%%%%%%%%%%%%%%%%%%%%%%%%%
%%%%%%%%%%%%%%%%%%%%%%%%%%%%%%%%%%%%%%%%%%%%%%%%%%

Corollary \ref{sl_2_5} $(ii)$, Theorems \ref{teor_3_2} $(i)$
and \ref{teor_4_6} imply  the following

%%%%%%%%%%%%%%%%%%%%%%%%%%%%%%%%%%%%%%%%%%%%%%%%%%
%%%%%%%%%%%%%%%%%%%%%%%%%%%%%%%%%%%%%%%%%%%%%%%%%%%%%%%%%%%%%%%%%%%%%%%%%%%%%%%%%%%%%%%%%%%%%%%%%%%%
\begin{theorems}\label{teor_4_7}
Let $\gamma$ and $\beta$ be oriented curves, generated by
non-degenerate $I$-path $x$ and $J$-path $y$ respectively, $u(s) =
x\big(q_x(s)\big)$, $s \in I(x)$, $v(r) = y\big(q_y(r)\big)$, $r \in
I(y)$. Then

(i) if the paths $u(s)$ and $v(r)$ are strongly regular, then the curves
$\gamma$ and $\beta$ are $O(n,p,K)$-equivalent $\big($respectively,
$SO(n,p,K)$-equivalent$\big)$, if and only if $I(x) = I(y)$ and
\begin{equation}\label{fla_22}
\left[\mathop{u}^{(m-1)}(s),\mathop{u}^{(m-1)}(s)\right] = \left[\mathop{v}^{(m-1)}(s),\mathop{v}^{(m-1)}(s)\right]
\end{equation}
$\big($respectively, equalities \eqref{fla_22} and
\[ \det M(u)(s) = \det M(v)(s)\big) \] are valid for all $s \in I(x)$ and $m = \overline{1,n}$ $\big($respectively,
$m= \overline{1,(n-1)}$ $\big)$;

(ii) if paths $u'(s)$ and $v'(s)$ are strongly regular, then curves
$\gamma$ and $\beta$ are $K^n \lhd O(n,p,K)$-equivalent
$\big($respectively, $K^n \lhd SO(n,p,K)$-equivalent$\big)$, if and
only if $I(x) = I(y)$ and equality \eqref{fla_22}
$\big($respectively, equalities \eqref{fla_22} and
\[ \det M'(u)(s) = \det M'(v)(s)\big) \] holds for all $s \in I(x)$, $m=\overline{2,(n+1)}$ $\big($respectively,
$m=\overline{2,n}$ $\big)$.
\end{theorems}
%%%%%%%%%%%%%%%%%%%%%%%%%%%%%%%%%%%%%%%%%%%%%%%%%%%%%%%%%%%%%%%%%%%%%%%%%%%%%%%%%%%%%%%%%%%%%%%%%%%%

\renewcommand{\refname}

\end{document}